\documentclass[12pt,english]{article}

\usepackage{anysize}
\marginsize{2cm}{2cm}{2.5cm}{2.4cm}

\usepackage[english]{babel}
\usepackage{amsfonts}
\usepackage{latexsym}
\usepackage{amsmath}
\usepackage{mathrsfs}
\usepackage{amssymb}

\let\phi=\varphi
\newcommand{\Z}{{\mathbb Z}}
\newcommand{\R}{{\mathbb R}}
\newcommand{\N}{{\mathbb N}}
\newcommand{\eps}{\varepsilon}
\newcommand{\IP}{{\mathbb P}}

\newcommand{\G}{{\mathcal G}}

\newcommand{\Po}{{\mathtt P}_{\omega}}
\newcommand{\Poi}{{\mathtt P}_{\omega'}}
\newcommand{\Eo}{{\mathtt E}_{\omega}}
\newcommand{\Eoi}{{\mathtt E}_{\omega'}}

\newcommand{\sig}{\sigma}
\newcommand{\gam}{\gamma}

\newcommand{\ga}{\text{\boldmath ${\gamma}$}}

\newcommand{\qed}{\hfill$\Box$\par\medskip\par\relax}
\newcommand{\1}[1]{{\mathbf 1}{\{#1\}}}

\newcommand{\dd}{{\mathtt d}}
\newcommand{\C}{\mathfrak{C}}

\let\phi=\varphi

\newcommand{\Bigexip}[1]{\Big\langle #1 \Big\rangle_{\!\!{}_\IP}}

\newcommand{\exip}[1]{\big\langle #1 \big\rangle_{\!{}_\IP}}

\newtheorem{theo}{Theorem}[section]
\newtheorem{lm}[theo]{Lemma}
\newtheorem{df}[theo]{Definition}
\newtheorem{prop}[theo]{Proposition}

\allowdisplaybreaks

\title{Random walks with unbounded jumps among random conductances I: Uniform quenched CLT}
\author{Christophe Gallesco \and Serguei
Popov}

\begin{document}

\bibliographystyle{plain}

\maketitle
{\footnotesize 

\noindent   Department of Statistics, 
Institute of Mathematics, Statistics and Scientific Computation,
University of Campinas -- UNICAMP, 
rua S\'ergio Buarque de Holanda 651, 13083--859, Campinas SP,
Brazil\\
\noindent e-mails: \texttt{gallesco@ime.unicamp.br}, \texttt{popov@ime.unicamp.br}

}

\maketitle

\begin{abstract}
We study a one-dimensional random walk among random conductances, 
with unbounded jumps. Assuming the ergodicity of the collection of
conductances and a few other technical conditions (uniform
ellipticity and polynomial bounds on the tails of the jumps)
we prove a quenched \textit{uniform} invariance principle for the
random walk. This means that 
the rescaled trajectory of length~$n$ is (in a certain sense)
close enough to the Brownian motion, uniformly with respect to
the choice of the starting location in an interval 
of length~$O(\sqrt{n})$ around the origin.
\\[.3cm] \textbf{Keywords:} ergodic environment, unbounded jumps, hitting
probabilities, exit distribution
\\[.3cm] \textbf{AMS 2000 subject classifications:} 60J10, 60K37
\end{abstract}

\section{Introduction and results}
\label{s_intro}
Suppose that for each pair of integers we are given a
nonnegative number. One may think that sites of~$\Z$ are nodes of an
electrical network where any site can be connected to any other 
site,
and those numbers are thought of as the \textit{conductances} of the
corresponding links. The conductances are initially chosen at
random, and we call the set of the conductances \textit{random
environment}. For the random environment, we assume that it is
stationary and ergodic.
 Given the conductances, one then defines a
(reversible) discrete-time random walk in the usual way: the transition
probability from~$x$ to~$y$ is proportional to the conductance
between~$x$ and~$y$. 

Here and in the companion paper~\cite{GP_2}
 we study one-dimensional random
walks among random conductances informally described above, with
unbounded jumps (we impose a condition that implies that 
the conductances can decay polynomially in the distance between the sites, but with sufficiently large power). 
The main result of~\cite{GP_2} 
concerns the (quenched) limiting law
of the trajectory of the random walk $(X_n, n=0,1,2,\ldots)$ starting
 from the origin up to time~$n$, under condition that it remains
positive at the moments $1,\ldots,n$. In~\cite{GP_2} we
prove that, after suitable rescaling, for a.e.\ environment it
converges to the
\emph{Brownian meander} process, which is, roughly speaking, a
Brownian motion conditioned on staying positive up to some finite
time. It turns out that one of the main ingredients 
for the proof of the conditional CLT is the \textit{uniform} quenched
CLT,
which is the main result (Theorem~\ref{Theouni}) of the present paper. 


Our main motivation for considering one-dimensional random walks with
unbounded jumps among random conductances and with minimal assumptions
on the environment
comes from Knudsen billiards in
random tubes, see \cite{CP11,CPSV1,CPSV2,CPSV3}.
 This model can be regarded as
a discrete-time Markov chain with continuous space, the positions
of the walker correspond to places where a billiard ball with a
\emph{random} law of reflection of certain form hits the boundary.
This model has some nice reversibility properties that makes it in
some sense a continuous space analogue of the random walk among
random conductances. In (rather long and technical) Section~3
of~\cite{CPSV3} the following problem was treated: 
given that the
particle, injected at the left boundary of the tube, crosses the
tube of length~$H$ without returning to the starting boundary, 
then the crossing time exceeds~$\eps H^2$ with high probability,
as $\eps\to 0$, $H\to\infty$.
Of course,
such a fact would be an easy consequence of a \emph{conditional} limit
theorem similar to the one described above
(since the probability that the Brownian meander
behaves that way is high). We decided to study
the discrete-space model because it presents less technical 
difficulties than random billiards in random tubes, and therefore allows 
to obtain finer results (such as the conditional CLT).

(Unconditional) quenched Central Limit Theorems for this and related
models (even in the many-dimensional case) received much attention
in the recent years, see e.g.~\cite{BD,BBK,Biskup,BB,BP,MP,M}.
Mainly, the modern approach consists in constructing a so-called
corrector function which turns the random walk to a martingale, and
then using the CLT for martingales. To construct the corrector, one
can use the method of orthogonal
projections~\cite{BP,MP,M}. While the corrector
method is powerful enough to yield quenched
CLTs, its construction by itself is not very explicit, and, in
particular, it does not say a lot about the speed of convergence.
Besides that, it is, in principle, not very clear how the speed of
convergence depends on the starting point. For instance, one may
imagine that there are ``distant'' regions where the environment is
``atypical'', and so is the behavior of the random walk starting
at a point from such a region until the time when it comes to
``normal''
parts of the environment. In Theorem~\ref{Theouni} we prove that,
for rather general ergodic environments that admit unbounded jumps,
if the rescaled trajectory by time~$n$ is ``close'' enough to the
Brownian motion, then so are the trajectories starting from points
of an interval of length~$O(\sqrt{n})$ centered in the origin. 
In our opinion, this
result is interesting by its own, but for us the main motivation for
investigating this question was that it provides an
important tool for proving the conditional CLT. Indeed,
the strategy of the proof of the conditional CLT is to force
it a bit away (around $\eps\sqrt{n}$) from the origin in a 
``controlled''
way and then use the usual (unconditional) CLT; but then it is clear 
that it is quite convenient to have the CLT for \emph{all}
starting positions in an interval of length~$O(\sqrt{n})$ at once.


It is important to note that in the most papers about random walks
with random conductances one assumes that the jumps are uniformly
bounded, usually nearest-neighbor (one can mention
e.g.~\cite{BBK,CFG} that consider the case of unbounded jumps). 
When there is no uniform bound on the size of the jumps, this of
course brings some complications to the proofs, as one can see in
the proof of Theorem~\ref{Theouni} below. Still, in our opinion, it
is important to be able to obtain ``usual'' results for the case of
long-range jumps as well; for example, in some related models, such
as the above-mentioned reversible random billiards in random 
tubes~\cite{CPSV2,CPSV3} the jumps are naturally unbounded.

In the case when (in dimension~$1$) the jumps are uniformly bounded,
the proofs become much simpler, mainly because one does not need to bother 
about the exit distributions, as in Section~\ref{s_estim_exit}.
The case of nearest-neighbor jumps is,
of course, even simpler, since many quantities of interest 
have explicit expressions. We will not discuss this case
separately, since it is (in some sense) ``too easy'' and
does not provide a lot of clues about how the walk with unbounded 
jumps should be treated. Let us make an observation that a random
walk with nearest-neighbor jumps becomes a very interesting and
complex object to study if one samples at random not the
conductances, but the transition probabilities themselves
(i.e., the transition probabilities from~$n$ to~$n+1$ are
chosen independently before the process starts). The resulting
random walk, while still reversible, behaves quite differently
(in particular, diffusive limits are unusual for that model). 
We only mention that \emph{conditional} (on being at the origin
at time~$2t$) behavior of this random walk in the transient 
case was studied in~\cite{GP09}, and a similar
result for the recurrent case can be obtained from
Corollary~2.1 of~\cite{CP}.

Of course, a natural question is whether the result analogous to
Theorem~\ref{Theouni} also holds for the many-dimensional 
nearest-neighbor random walk among random conductances.
We postpone the discussion about that to the end of this section.

Now, we define the model formally.
For $x,y \in \Z$, let us denote by $\omega_{x,y}=\omega_{y,x}$ the
conductance
between~$x$ and~$y$. Define $\theta_z\omega_{x,y}:=\omega_{x+z,y+z}$, for all $z\in \Z$.
Note that, by Condition~K below, the vectors~$\omega_{x,\cdot}$ are
elements of the Polish space $\ell^2(\Z)$.
We assume that $(\omega_{x,\cdot})_{x\in \Z}$ is a stationary 
ergodic (with respect to the family of shifts $\theta$) sequence of
random vectors; $\IP$ stands for the law of this sequence
and~$\exip{\cdot}$ is the
expectation with respect to~$\IP$. The collection of all
conductances $\omega=(\omega_{x,y}, x, y\in \Z)$ is called the
\textit{environment}. For all $x\in\Z$, define
$C_x:=\sum_{y}\omega_{x,y}$. Given that $C_x<\infty$ for all $x\in\Z$
(which is always so by Condition~K below),
the random walk~$X$ in random environment~$\omega$ is defined through
its transition probabilities
\[
p_{\omega}(x,y)=\frac{\omega_{x,y}}{C_x};
\]
that is, if~$\Po^x$ is the quenched law of the random walk starting
from $x$, we have 
\[
 \Po^x[X_0=x]=1, \quad \Po^x[X_{k+1}=z\mid X_k=y]=p_{\omega}(y,z).
\]
Clearly, this random walk is reversible with the reversible
measure~$(C_x,x\in\Z)$.
 Also, we denote
by~$\Eo^x$ the quenched expectation for the process starting
 from~$x$. When the random walk starts from~$0$, we use shortened
notations $\Po,\Eo$.

In order to prove our results, we need to make two
technical assumptions on the environment:

\medskip
\noindent
{\textbf {Condition~E}.}
There exists $\kappa>0$ such that, $\IP$-a.s., 
$\omega_{0,1}\geq \kappa$.

\medskip
\noindent{\textbf{Condition~K}.} 
There exist constants $K,\beta>0$ such that, $\IP$-a.s.,
$\omega_{0,y}\leq \frac{K}{1+y^{3+\beta}}$,
for all $y\geq 0$.
Note, for future reference, that the stationarity of $\IP$ and Conditions E and K together imply that there exists $\hat{\kappa}>0$ such that, $\IP$-a.s., 
\begin{equation}
\label{Ellip2}
\hat{\kappa} \leq \sum_{y\in \Z} \omega_{0,y}\leq 
\hat{\kappa}^{-1}.
\end{equation}

\medskip

We decided to 
formulate Condition~E this way because, due to the fact that
this work was motivated by random billiards, the main challenge was
to deal with the long-range jumps.
It is plausible that Condition~E could be relaxed to some extent; 
however, for the sake of cleaner presentation of the argument, 
we prefer not trying to deal with \emph{both} long-range jumps and
the lack of nearest-neighbor ellipticity.


Next, for all $n\geq 1$, we define the continuous map $Z^n=(Z^n(t),t\in 
\R_+)$ as the natural polygonal interpolation of the map $k/n\mapsto
\sig^{-1}n^{-1/2}X_k$ (with $\sig$ from Theorem~\ref{t_q_invar_princ}
below). In other words,
\[
\sig \sqrt{n}Z^n_t = X_{\lfloor nt\rfloor}+(nt-\lfloor
nt\rfloor)X_{\lfloor nt\rfloor+1}
\]
with $\lfloor \cdot\rfloor$ the integer part. Also, we denote by~$W$
the standard Brownian motion.

First, we state the following result, which is the usual quenched
invariance principle:
\begin{theo}
\label{t_q_invar_princ}
Assume Conditions~E and~K.
 Then, there exists a finite (nonrandom) constant $\sigma>0$ such
that for $\IP$-almost all~$\omega$, $Z^n$
converges in law, under $\Po$, to Brownian motion~$W$ as
$n\to\infty$.
\end{theo}
Of course, with the current state of the art in this
field, obtaining the proof of Theorem~\ref{t_q_invar_princ} is a
mere exercise (one can follow e.g.\ the argument of~\cite{BP});
for this reason, we do not write the proof of
Theorem~\ref{t_q_invar_princ} in this paper. The
key observation, though, is that Condition~K implies that 
\[
 \Bigexip{\sum_{y\in \Z} y^2\omega_{0,y}} < \infty.
\]

Let $C(\R_+)$ be the space of continuous functions
from $\R_+$ into $\R$. Let us denote by $\C_b(C(\R_+),\R)$ (respectively,
$\C^u_b(C(\R_+),\R)$) the space of bounded
continuous (respectively, bounded uniformly continuous)
functionals from $C(\R_+)$ into $\R$ and by $\mathcal{B}$ the Borel
$\sig$-field on $C(\R_+)$. We have the following result, which is
referred to as quenched 
Uniform Central Limit Theorem (UCLT):
\begin{theo}
\label{Theouni}
Under Conditions~E and~K, the following statements hold:
\begin{itemize}
\item[(i)] we have $\IP$-a.s., for all $H>0$ and any 
$F\in \C_b(C(\R_+),\R)$,
\[
\lim_{n \to \infty} \sup_{x\in [-H\sqrt{n},H\sqrt{n}]}\Big| 
 {\mathtt E}_{\theta_x \omega}[F(Z^{n})]-E[F(W)]\Big|=0;
\]
\item[(ii)] we have $\IP$-a.s., for all $H>0$ and any
$F\in \C^u_b(C(\R_+),\R)$, 
\[
\lim_{n \to \infty} \sup_{x\in [-H\sqrt{n},H\sqrt{n}]}\Big| 
 {\mathtt E}_{\theta_x \omega}[F(Z^{n})]-E[F(W)]\Big|=0;
 \]
 \item[(iii)] we have $\IP$-a.s., for all $H>0$ and any closed 
set~$B\in\mathcal{B}$,
\[
\limsup_{n \to \infty} \sup_{x\in [-H\sqrt{n},H\sqrt{n}]}
 {\mathtt P}_{\theta_x \omega}[Z^{n}\in B]\leq P[W\in B];
 \]
\item[(iv)] we have $\IP$-a.s., for all $H>0$ and any open set~$G\in\mathcal{B}$, 
\[
\liminf_{n \to \infty} \inf_{x\in [-H\sqrt{n},H\sqrt{n}]}
{\mathtt P}_{\theta_x \omega}[Z^{n}\in G]\geq P[W\in G];
 \]
 \item[(v)] we have $\IP$-a.s., for all $H>0$ and any $A\in
\mathcal{B}$ such that $P[W\in \partial A]=0$,
\[
\lim_{n \to \infty} \sup_{x\in [-H\sqrt{n},H\sqrt{n}]}\Big| 
 {\mathtt P}_{\theta_x \omega}[Z^{n}\in A]-P[W\in A]\Big|=0.
 \]
\end{itemize}
\end{theo}

Even though it may be possible to find a concise formulation
of our main result with only one ``final'' statement
and not a list of equivalent ones (the authors did not succeed in finding it), 
we content ourselves in writing it in this form
because, in our opinion, possible situations where it can be
useful are covered by the list.
Of course, item~(ii) is redundant (it follows trivially from~(i)),
and (iii) and (iv) are equivalent by complementation.
 
For the model of this paper, it is not possible to
generalize Theorem~\ref{Theouni} by considering a wider interval
$[-Hn^\alpha,Hn^\alpha]$ for some~$\alpha>1/2$,
where the starting point is taken; this is
because we only assume the
ergodicity of the environment of conductances. Indeed, 
consider e.g.\ a nearest-neighbor
random walk, and suppose that the conductances can assume only two possible 
values, say, $1$ and~$2$. To construct the stationary ergodic 
random environment, we first construct its cycle-stationary version 
in the following way. Fix~$\eps>0$ such that $\frac{\alpha}{1+\eps}>\frac{1}{2}$,
and let us divide the edges of~$\Z$ into blocks of random i.i.d.\ sizes 
$(V_i, i\in\Z)$, with $P[V_1>s]=O(s^{-(1+\eps)})$. Inside each block, 
we toss a fair coin and, depending on the result,
 either place all $1$s, or an alternating sequence of $2$s and $1$s. 
Since the expected size of the block is finite, it is clear that this 
environment
can be made stationary (and, of course, ergodic) by a standard
random shift procedure, see e.g.\ Chapter~8 of~\cite{Th}.
Then, one readily obtains that in the interval $[-Hn^\alpha,Hn^\alpha]$,
with large probability, one finds both $1\ldots 1$-blocks
and $212\ldots12$-blocks of length at least $\sqrt{n}$. So, the UCLT
cannot be valid: just consider starting points in the middle of 
two blocks of different type. 
 It is, in our opinion, an interesting problem to obtain a
stronger form of Theorem~\ref{Theouni} in the case
when the environment has mixing or independence properties.
It seems plausible that one can make the above interval at least
polynomially (with any power) wide, but we prefer not to discuss
further questions of this type in this paper:
 in any case, for the results of~\cite{GP_2},
 Theorem~\ref{Theouni} is already enough.

Let us also comment on possible many-dimensional
variants of Theorem~\ref{Theouni}. For the case of nearest-neighbor
random walks in~$\Z^d$ with random conductances bounded from both sides
by two positive constants, an analogous result was obtained
in~\cite{GGPV} (Theorem~1.1). The proof of Theorem~1.1 of~\cite{GGPV}
relies on the uniform heat-kernel bounds of~\cite{Del}; one
uses these bounds to obtain that, regardless of the starting point,
with probability close to~$1$ the walk will enter to the set of ``good''
sites (i.e., the sites from where the convergence is good enough). 
Naturally,
this poses the question of what to do with unbounded conductances
(and/or unbounded jumps),
to which we have no answer for now
(although one can expect, as usual, that the case $d=2$ should be more
accessible, since in this case each site is ``surrounded'' by ``good''
sites, cf.\ e.g.\ the proof of Theorem~4.7 in~\cite{Biskup}).

The paper is organized in the following way:
in the next section, we obtain some auxiliary facts which 
are necessary for the proof of Theorem~\ref{Theouni}
(recurrence, estimates on the probability of confinement in an interval,
estimates on the exit measure from an interval). Then, in
Section~\ref{s_proof_t_uni}, we give the
proof of Theorem~\ref{Theouni}.

We will denote by $K_1$, $K_2$,~$\dots$ the ``global'' 
constants, that is, those that are used all along the paper, and by
$\gam_1$, $\gam_2$,~$\dots$ the ``local" constants, that is, those
that are used only in the subsection in which they appear for the
first time. For the local constants, we restart the numeration in the
beginning of each subsection. 
Depending on the context, expressions like 
$x\in [-H\sqrt{n},H\sqrt{n}]$ should be understood as $x\in
[-H\sqrt{n},H\sqrt{n}]\cap \Z$.

\section{Auxiliary results}
\label{s_aux_results}
In this section, we will prove some technical results that 
will be needed later to prove Theorem~\ref{Theouni}.
Let us introduce the following notations. If $A\subset \Z$,
\[
\tau_A:=\inf\{n\geq 0: X_n\in A\} \quad \text{and} \quad
\tau^+_A:=\inf\{n\geq 1: X_n\in A\}.
\]

\subsection{Recurrence of the random walk}
\label{s_recurr_RW}
\begin{lm}
\label{Rec}
Under Conditions E and K the random walk $X$ is $\IP$-a.s.\ 
recurrent.
\end{lm}
\textit{Proof.} To show the recurrence of the random walk,
 we will show that the probability of escape to infinity is zero.
First, let us consider the finite interval $I_L=[-L,L]$ for some
$L>0$. Consider the time~$\tau_{I^c_L}$ of exit from the 
interval~$I_L$.
By the Dirichlet variational principle for reversible 
Markov chains (see, for example, Theorem 6.1 Chap.~II of~\cite{Lig}) we have
that 
\begin{equation}
\label{JFK}
2C_0\Po[\tau_{I^c_L}<\tau^+_0]=\min_{f\in \mathcal{H}}\Phi(f)
\end{equation}
where $\Phi$ is the Dirichlet form defined by
\[
\Phi(f):=\sum_{x,y\in \Z}\omega_{x,y}[f(x)-f(y)]^2
\]
and $\mathcal{H}$ is the following set of functions:
\[
\mathcal{H}:=\{f:\Z\to [0,1]: f(0)=0
\text{ and } f(x)=1  \text{ for }
x\notin I_L\}.
\]
In order to estimate $\Po[\tau_{I^c_L}<\tau^+_0]$ from 
above let us consider the function $h$ in $\mathcal{H}$
defined by 
\[
h(x)= \left\{
    \begin{array}{ll}
       L^{-1}|x|, & \mbox{if}~|x|\leq L,\\
       1, & \mbox{if}~|x|>L.\\
    \end{array}
\right.
\]
Now, let us estimate $\Phi(h)$.
We start by writing
\begin{align}
\label{Dirich0}
\Phi(h)
&=\sum_{x,y\in \Z}\omega_{x,y}[h(x)-h(y)]^2\nonumber\\
&=\sum_{-L< x,y < L}\omega_{x,y}[h(x)-h(y)]^2+2\sum_{x\in (-\infty, -L]\cup [L,\infty)}\sum_{ y\in (-L,L)}
 \omega_{x,y}[h(x)-h(y)]^2.
\end{align}
We are going to show that both terms in the decomposition (\ref{Dirich0}) of $\Phi(h)$ are of order smaller or equal to~$L^{-1}$.
Indeed, it is not difficult to see that each of both terms in the decomposition (\ref{Dirich0}) is smaller than 
\begin{equation*}
\frac{2}{L^2}\sum_{-L<x<L}\sum_{y\in \Z} 
\omega_{x,y}(y-x)^2.
\end{equation*}
By Condition~K, we obtain that 
there exists a constant $\gam_1$ such that $\IP$-a.s., $\sum_{y\in
\Z} \omega_{x,y}(y-x)^2\leq\gam_1$ for all~$x$.
Then, we deduce
\begin{equation*}
\Phi(h)\leq 
\frac{8\gam_1}{L}.
\end{equation*}
Using (\ref{JFK}), we obtain that,
\begin{equation*}
2C_0\Po[\tau_{I^c_L}<\tau^+_0]\leq \frac{8 \gamma_1}{L}.
\end{equation*}
By (\ref{Ellip2}), we have $C_0\geq \hat{\kappa}$ so that
\begin{equation}
\label{decro}
\Po[\tau_{I^c_L}<\tau^+_0]\leq \frac{4\gamma_1}{\hat{\kappa}L}.
\end{equation}
Now, let~$p_{esc}$ be the probability that the walk started 
at~$0$ escapes to infinity. We have $p_{esc}=\lim_{L\to
\infty}\Po[\tau_{I^c_L}<\tau^+_0]$, and so, 
 by~(\ref{decro}), $p_{esc}=0$.
Hence, the random walk~$X$ is $\IP$-a.s.\ recurrent.
\qed

\subsection{Probability of confinement}
\label{s_prob_conf}
Let $I=[a,b]\subset \Z$ be a finite interval containing at least 
3 points and let $B=(-\infty, a]$ and $E=[b,\infty)$. In this
subsection we shall prove the following
\begin{prop}
\label{specgap}
There exist constants $K_1>0$ and $K_2>0$ such that we have
$\IP$-a.s.,
\begin{equation*}
\max_{x\in (a,b)}\Po^x[\tau_{B\cup E}> n]\leq \exp 
\Big\{-\frac{n}{K_1(b-a)^2}\Big\}
\end{equation*}
for all $n> K_2(b-a)^2$.
\end{prop}
\textit{Proof.} Let $\omega$ be a realization of the 
random environment. Consider the new environment obtained from
$\omega$ by deleting all the conductances  $\omega_{x,y}$ if~$x$ 
and~$y$ belong to $B\cup E$. The reversible measure $C'$ on this new
environment is given by
\begin{align*}
C'_x=\left\{
    \begin{array}{ll}
       C_x, &\mbox{if $x\in (a, b)$,}\\
       \sum_{y\notin B\cup E}\omega_{x,y},& \mbox{otherwise}.\\
    \end{array}
\right.
\end{align*}
Next, we define $C'_B:=\sum_{y\in B}C'_y$ and $\pi_B(x):=C'_x/C'_B$ 
for all $x \in B$.
Observe that, by Conditions~E and K, $C'_B$ is positive and 
finite $\IP$-a.s. Hence, it holds that~$\pi_B$ is $\IP$-a.s.\ a
probability measure on~$B$.
In the same way we define the probability measure~$\pi_E$ on~$E$.
Now, we introduce a new Markov chain $X'$ on 
a finite state space $\mathcal{S}':=(a,b)\cup\{\Delta_B,
\Delta_E\}$ ($\Delta_B$ and $\Delta_E$ are the states
corresponding to~$B$ and~$E$). 
On $\mathcal{S}'$, we define the following transition
probabilities:
if $x\notin \{\Delta_B,\Delta_E\}$ 
\begin{equation*}
P_{x,\Delta_E}=\sum_{y\in E}\frac{\omega_{x,y}}{C'_x},
\phantom{**}P_{x,\Delta_B}=\sum_{y\in B}\frac{\omega_{x,y}}{C'_x}
\end{equation*} 
and
\begin{equation*}
P_{\Delta_E,x}=\sum_{y\in E}\pi_E(y)\frac{\omega_{x,y}}{C'_y},
\phantom{**}P_{\Delta_B,x}=\sum_{y\in
B}\pi_B(y)\frac{\omega_{y,x}}{C'_y}.
\end{equation*} 
Then, set
$P_{\Delta_E,\Delta_B}=P_{\Delta_B,\Delta_E}
 =P_{\Delta_B,\Delta_B}=P_{\Delta_E,\Delta_E}=0$.
For $x \notin \{\Delta_B,\Delta_E\}$ and 
$y\notin \{\Delta_B,\Delta_E\}$ we just set
$P_{x,y}=\omega_{x,y}/C'_x$.
Defining $C'_{\Delta_B}:=C'_B$ and $C'_{\Delta_E}:=C'_E$, 
we can easily check that the detailed balance equations are
satisfied, that is, on $\mathcal{S}'$ we have a new set of
conductances~$\omega'$ specified by
$\omega'_{x,y}:=C'_xP_{x,y}=C'_yP_{y,x}$. Observe also that by Condition~K, there exists 
a constant $\gam_1>0$ such that $\IP$-a.s., $C'_x\leq \gam_1$
for all $x\in \mathcal{S}'$. 
By the commute time identity 
(see for example Proposition 10.6 of~\cite{LevPeres}) we have that 
\[
\Eo^x[\tau_{B\cup E}]\leq \Eoi^x[\tau_{\Delta_{B}}]+ \Eoi^{\Delta_{B}}[\tau_{x}]=\Big( \sum_{y\in \mathcal{S}'} C'_y\Big) R_{\text{eff}}(\Delta_B, x)
\]
where $R_{\text{eff}}(\Delta_B, x)$ is the effective resistance between $\Delta_B$ and $x$. We have 
\begin{equation*}
\Big( \sum_{y\in \mathcal{S}'} C'_y\Big) \leq \gam_1(b-a+1)
\end{equation*}
and 
\begin{equation*}
R_{\text{eff}}(\Delta_B, x)\leq \sum_{y=\Delta_B}^{x-1}\omega^{-1}_{y,y+1} \leq \kappa^{-1}(b-a+1).
\end{equation*}
Thus,
\begin{equation*}
\Eo^x[\tau_{B\cup E}]\leq \gam_1\kappa^{-1}(b-a+1)^2\leq \gam_2(b-a)^2
\end{equation*}
for some positive constant $\gam_2$.
By the Chebyshev inequality, we can choose a large enough constant $\gam_3>0$ in such a way that
\begin{equation}
\label{rtyu}
\Po^x\big[\tau_{B\cup E}> \lfloor\gam_3(b-a)^2\rfloor\big]\leq \frac{\Eo^x[\tau_{B\cup E}]}{\lfloor \gam_3 (b-a)^2\rfloor}\leq \frac{\gam_2(b-a)^2}{\lfloor \gam_3 (b-a)^2\rfloor} <1.
\end{equation}
Let us denote $s:=\lfloor\gam_3(b-a)^2\rfloor$ and $p:= \gam_2(b-a)^2 \lfloor \gam_3 (b-a)^2\rfloor^{-1}$. For $n\geq s$ divide the time interval $[0,n]$ into 
$N:=\lfloor\frac{n}{s}\rfloor$ subintervals of length~$s$.
Using~(\ref{rtyu}) and the Markov property we obtain
\begin{align*}
\label{conttimeestim}
 \Po^x[\tau_{B\cup E}> n]
&\leq \Po^x\big[X'(sj)\notin \{\Delta_B, \Delta_E\}, j=1,\dots,N\big] \nonumber\\
&\leq (1- p)^N\nonumber\\
&\leq \exp \Big(-\frac{n}{\gam_{4}(b-a)^2}\Big)
\end{align*}
for some positive constant $\gam_4$.
This concludes the proof of Proposition~\ref{specgap}. 
\qed

\subsection{Estimates on the exit distribution}
\label{s_estim_exit}
Let $I=[a,b]\subset \Z$ be a finite interval and $E=(-\infty, a]\cup[b,+\infty)$. We prove the following
\begin{prop}
\label{SNL}
For all $\eta>0$ there exists $M>0$ such that $\IP$-a.s., 
for each interval $[a,b]\subset \Z$ containing at least three points 
we have
\[
\min_{x\in (a,b)}\Po^x[X_{\tau_{E}}\in I_M]\geq 1-\eta
\]
with $I_M:=[a-M,a]\cup [b,b+M]$.
\end{prop}
\textit{Proof.} Fix an arbitrary $\eta\in (0,1)$.  
For intervals $[a,b]$ of length 2, there exists only one point $x$ in
$(a,b)$. By the Markov property we have that
\begin{equation*}
\Po^x[X_{\tau_{E}}\in I_{M}]= \Po^x[X_1\in I_M\mid X_1\neq x].
\end{equation*}
This implies that
\begin{equation*}
\Po^x[X_{\tau_{E}}\in I_{M}]= 1-\frac{\Po^x[X_1\in (-\infty, a-M)\cup(b+M,\infty)]}{\Po^x[X_1\neq x]}.
\end{equation*}
Then, Condition K and (\ref{Ellip2}) guarantee the existence of a constant~$M>0$ such
that $\IP$-a.s.,
\begin{equation*}
\Po^x[X_{\tau_{E}}\in I_{M}]\geq 1-\eta.
\end{equation*}

For intervals $[a,b]$ of length greater or equal to 3, 
let us do the following. 
Fix some
$x\in (a,b)$. Let $\zeta_0=0$ and for $i\geq 1$, $\zeta_i:=\inf\{n>\zeta_{i-1}: X_n=x\}$
with the convention $\inf\{\emptyset\}=+\infty$.
Since by Lemma~\ref{Rec} our random walk is $\IP$-a.s.\ recurrent,
the sequence $(\zeta_i)_{i\geq1}$ is $\IP$-a.s.\ strictly increasing
and we have by the Markov property
\begin{align}
\label{Simplif}
\Po^x[X_{\tau_E}\in I_M]
&=\sum_{i=0}^{\infty}\Po^x[X_{\tau_E}\in I_M \mid \tau_{E} \in 
[\zeta_i, \zeta_{i+1})]\Po^x[\tau_{E} \in [\zeta_i,
\zeta_{i+1})]\nonumber\\
&=\Po^x[X_{\tau_E}\in I_M \mid \tau_{E}<\tau_x^+].
\end{align}
Let us define $A_{E}:=\{\tau_E<\tau_x^+\}$.
First, we write
\begin{align}
\label{XM1}
\Po^x[X_{\tau_{E}}\in I_M\mid A_E]
&=1-\Po^x[X_{\tau_{E}}\notin I_M\mid A_E]\nonumber\\
&=1-\sum_{y\in(-\infty,a-M)\cup(b+M,\infty)} 
\Po^x[X_{\tau_{E}}=y\mid A_E].
\end{align}
Then, consider the new environment~$\omega'$ obtained from~$\omega$
by deleting all the conductances $\omega_{y,z}$ when both~$y$ and~$z$
belong to~$E$. The reversible measure on this new 
environment~$\omega'$ is given by
\begin{align*}
C'_y:= \left\{
    \begin{array}{ll}
    C_y, & \mbox{if $y\in (a, b)$,}\\
    \sum_{z\notin E}\omega_{y,z}, & \mbox{otherwise}.\\
    \end{array}
 \right. 
\end{align*}
We define $C'_E:=\sum_{y\in E}C'_y$ and for all $y\in E$, $\pi_E(y):=C'_y/C'_E$.
Observe that by Condition~K, 
$C'_E\in (0,\infty)$, $\IP$-a.s. Hence, $\pi_E$ is a probability
measure on~$E$. For the sake of simplicity we write~$\Poi^E$
instead of $\Poi^{\pi_E}$ for the random walk on $\omega'$ starting
with probability $\pi_E$.
We can couple the random walks in the environments~$\omega$ and
$\omega'$ so that $\Poi^x[X_{\tau_{E}}=y\mid
A_E]=\Po^x[X_{\tau_{E}}=y\mid A_E]$. 

Now, let us find an upper bound for the term 
$\Poi^x[X_{\tau_{E}}=y\mid A_E]$ with
$y\in(-\infty,a-M)\cup(b+M,\infty)$.
By definition of $A_E$ we have
\begin{equation}
\label{XM2}
\Poi^x[X_{\tau_{E}}=y \mid A_E]=\frac{\Poi^x[X_{\tau_{E}}=y,
\tau_E<\tau_x^+]}{\Poi^x[\tau_E<\tau_x^+]}.
\end{equation}
Let us denote by $\Gamma_{z',z''}$ the set of finite paths 
$(z',z_1,\ldots,z_k,z'')$ such that $z_i\notin E\cup\{z',z''\}$
for all $i=1,\ldots,k$.
Let $\ga=(z',z_1,\dots,z_k,z'')\in \Gamma_{z',z''}$ and
define 
\[
 \Poi^{z'}[\ga]:=\Poi^{z'}[X_1=z_1,\dots,X_k=z_k,X_{k+1}=z''].
\]
By reversibility we obtain
\begin{equation*}
\Poi^x[X_{\tau_{E}}=y,\tau_E<\tau_x^+]
=\sum_{\ga \in \Gamma_{x,y}}\Poi^x[\ga] 
=\frac{1}{C'_x}\sum_{\ga \in \Gamma_{y,x}}C'_y\Poi^y[\ga]
=\frac{C'_y}{C'_x}\Poi^y[\tau_x<\tau_E^+]
\end{equation*}
and
\begin{align}
\label{Relblar}
\Poi^x[\tau_E<\tau_x^+]
&=\sum_{z\in E}\sum_{\ga\in \Gamma_{x,z}}\Poi^x[\ga] 
=\sum_{z\in E} \sum_{\ga\in \Gamma_{z,x}}\frac{C'_z}{C'_x}
\Poi^z[\ga] 
=\frac{C'_E}{C'_x}\sum_{z\in E} \pi_E(z)
\sum_{\ga\in \Gamma_{z,x}}\Poi^z[\ga]\nonumber\\
&=\frac{C'_E}{C'_x}\Poi^{E}[\tau_x<\tau_E^+].
\end{align}
Thus, by~(\ref{XM2}) we have
\begin{equation}
\label{VL0}
\Poi^x[X_{\tau_{E}}=y \mid A_E]
=\frac{C'_y\Poi^y[\tau_x<\tau_E^+]}{C'_E\Poi^{E}[\tau_x<\tau_E^+]}.
\end{equation}
To bound from below the term $\Poi^{E}[\tau_x<\tau_E^+]$ we 
use an electric networks argument. To this end, we will define a
Markov chain on a new state space for which it will be easy to
compute the effective conductance. First, we introduce a point
$\Delta_E$ and the state space $\mathcal{S}:=(\Z \setminus
E)\cup\{\Delta_E\}$. For $z\notin E$, we define the transition
probabilities 
\begin{equation*}
P_{z,\Delta_E}=\sum_{u\in E}\frac{\omega'_{z,u}}{C'_z},
\qquad
P_{\Delta_E,z}=\sum_{u\in E}\pi_E(u)\frac{\omega'_{z,u}}{C'_u}.
\end{equation*} 
For $z\notin E$ and $u\notin E$ we set
$P_{z,u}=\omega'_{z,u}/C'_z$,
and, we put
$P_{\Delta_E,\Delta_E}=0$.
By defining $C'_{\Delta_E}:=C'_E$, we can easily check that 
the detailed balance equations are satisfied, i.e., for all $z\in
\mathcal{S}$ we have $C'_zP_{z,u}=C'_uP_{u,z}$.
We have that 
\begin{equation}
\label{Condeff1}
\Poi^{E}[\tau_{x}<\tau^+_E]
 =\Poi^{\Delta_E}[\tau_{x}<\tau^+_{\Delta_E}]=\frac{C_{\text{eff}}
(\Delta_E,x)}{C'_{\Delta_E}}
 =\frac{C_{\text{eff}}(\Delta_E,x)}{C'_{E}
}
\end{equation}
where $C_{\text{eff}}(\Delta_E,x)$ is the 
effective conductance between $\Delta_E$ and $x$.
Observe that
\begin{equation*}
C_{\text{eff}}(\Delta_E,x)\geq 
\Big(\sum_{i=a}^{x-1} \omega_{i,i+1}^{-1}\Big)^{-1} +
 \Big(\sum_{i=x}^{
b-1}\omega_{i,i+1}^{-1}\Big)^{-1}.
\end{equation*}
Using Condition~E, we obtain
\begin{equation}
\label{WQES}
C'_E\Poi^{E}[\tau_x<\tau_E^+]\geq \kappa
 \Big(\frac{1}{x-a}+\frac{1}{b-x}\Big).
\end{equation}
Then, we have to treat the term $C'_y\Poi^y[\tau_x<\tau_E^+]$.  
By construction of~$\omega'$
\begin{equation*}
C'_y\Poi^y[\tau_x<\tau_E^+]
=C'_y\sum_{z\in (a,b)}\!p_{\omega'}(y,z)\Poi^z[\tau_x<\tau_E]
=\sum_{z\in (a,b)}\!\omega'_{y,z}\Poi^z[\tau_x<\tau_E]
=\sum_{z\in (a,b)}\!\omega_{y,z}\Po^z[\tau_x<\tau_E].
\end{equation*}
Finally, we have to estimate $\Po^z[\tau_x<\tau_E]$ 
for $z\in (a,b)\setminus \{x\}$.
To this end, we define the following sequence of stopping times. 
Let $\Upsilon_0=0$ and for $i\geq 1$, $\Upsilon_i:=\inf\{n>\Upsilon_{i-1}: X_n=z\}$
with the convention $\inf\{\emptyset\}=+\infty$.
The sequence $(\Upsilon_i)_{i\geq1}$ is
a.s.\ strictly increasing and we have
\begin{align}
\label{VL1}
\Po^z[\tau_x<\tau_E]
=\Po^z[\tau_x<\tau_E 
 \mid \tau_{E\cup\{x\}} \in [0, \Upsilon_{1})].
\end{align}
Then, we have
\begin{equation}
\label{VL2}
\Po^z[\tau_x<\tau_E \mid \tau_{E\cup\{x\}} \in [0, \Upsilon_{1})]
=\frac{\Po^z[\tau_x<\tau_E, \tau_{E\cup\{x\}} 
\in [0, \Upsilon_{1})]}{\Po^z[\tau_{E\cup\{x\}} \in [0,
\Upsilon_{1})]}
\leq \frac{\Po^z[\tau_x<\tau_z^+]}{\Po^z[\tau_E<\tau_z^+]}\wedge 1.
\end{equation}
We estimate $\Po^z[\tau_E<\tau_z^+]$ 
in the following way,
\begin{align}
\label{VL3}
\Po^z[\tau_E<\tau_z^+]
&=\frac{C_{\text{eff}}(z,E)}{C_z}\nonumber\\
&\geq \frac{1}{C_z}
\Big(\Big(\sum_{i=z}^{b-1}\omega_{i,i+1}^{-1}\Big)^{-1} 
+\Big(\sum_{i=a}^{z-1}\omega_{i,i+1}^{-1}\Big)^{-1}\Big)
\nonumber\\
&\geq \hat{\kappa} \kappa\Big(\frac{1}{b-z}+\frac{1}{z-a} \Big).
\end{align}
Now, using the Dirichlet variational principle, we obtain an upper 
bound for $\Po^z[\tau_x<\tau_z^+]$.
Suppose that the interval $(x, b)\neq \emptyset$ and that $z\in
(x,b)$, consider the function~$h$ given by 
\[
h(u)= \left\{
    \begin{array}{ll}
       1, & \mbox{if}~u<x\phantom{*}\mbox{or}\phantom{*}u>2z-x,\\
       \frac{|z-u|}{z-x}, & \mbox{if}~x\leq u\leq z.
    \end{array}
\right.
\]
Hence, we have $2C_z\Po^z[\tau_x<\tau_z^+]\leq \Phi(h)$. 
By the same reasoning as we used in order to obtain (\ref{decro}) in the proof of Lemma \ref{Rec}, we deduce that there exists a constant~$\gam_2>0$ such that
\begin{equation}
\label{VL4}
\Po^z[\tau_x<\tau_z^+]\leq \frac{\gam_2}{z-x}
\end{equation}
for $z\in (x,b)$. Similarly, if we suppose that 
$(a,x)\neq \emptyset$ and $z\in (a,x)$, we obtain a bound similar
to~(\ref{VL4}) for $\Po^z[\tau_x<\tau_z^+]$. Then, we obtain that for
$z\in (a,b)\setminus \{x\}$,
\begin{equation}
\label{VLH}
\Po^z[\tau_x<\tau_z^+]\leq \frac{\gam_2}{|z-x|}.
\end{equation}
Note that we can choose~$\gam_2$ in such a way that it does not
depend on the size of the interval~$[a,b]$.
By~(\ref{VL1}), (\ref{VL2}), (\ref{VL3}) and~(\ref{VLH}) we
obtain
\begin{equation*}
C'_y\Poi^y[\tau_x<\tau_E^+]\leq
\sum_{ z\in (a,b)}\omega_{y,z}
 \Big(\frac{\gam_2(z-a)(b-z)}{\hat{\kappa} \kappa
(b-a)|z-x|}\wedge1\Big).
\end{equation*}
Thus, by~(\ref{VL0}) and (\ref{WQES}) we obtain
\begin{align}
\Poi^x[X_{\tau_{E}}=y \mid A_E] &\leq 
 \frac{1}{ \kappa}\frac{(x-a)(b-x)}{b-a}\sum_{ z\in
(a,b)}\omega_{y,z}\Big(\frac{\gam_2(z-a)(b-z)}{\hat{\kappa} \kappa
(b-a)|z-x|}\wedge1\Big), \nonumber\\
\Po^x[X_{\tau_{E}}\in I_M\mid A_E]
&\geq 1-\frac{1}{\kappa}\sum_{y\in I'_M}
 \frac{(x-a)(b-x)}{b-a}\sum_{z\in
(a,b)}\omega_{y,z}\Big(\frac{\gam_2(z-a)(b-z)}{\hat{\kappa} \kappa
(b-a)|z-x|}\wedge1\Big)\label{Ertac}
\end{align}
with $I'_M=(-\infty, a-M)\cup(b+M,\infty)$.
Let us divide the set~$I'_M$ into the subintervals
$J_1(M)=(b+M,\infty)$ and $J_2(M)=(-\infty, a-M)$. Denote
\begin{align*}
H_1(M)&=\sum_{y\in J_1}\frac{(x-a)(b-x)}{b-a}
 \sum_{ z\in
(a,b)}\omega_{y,z}
\Big(\frac{\gam_2(z-a)(b-z)}{\hat{\kappa} \kappa
(b-a)|z-x|}\wedge1\Big),\\
H_2(M)&=\sum_{y\in J_2}\frac{(x-a)(b-x)}{b-a}\sum_{ z\in (a,b)}
\omega_{y,z}\Big(\frac{\gam_2(z-a)(b-z)}{\hat{\kappa} \kappa
(b-a)|z-x|}\wedge1\Big).
\end{align*}
We have
\begin{align*}
H_1(M)
&\leq \sum_{y\in J_1}\frac{(x-a)(b-x)}{b-a}
\Big\{\sum_{ z\in
(a,\frac{x+a}{2}]}\frac{\gam_2(z-a)(b-z)}{\hat{\kappa} \kappa
(b-a)|z-x|} \omega_{y,z}
+\sum_{ z\in (\frac{x+a}{2},\frac{b+x}{2})} 
\omega_{y,z}
\nonumber\\
&\phantom{*****************}
+\sum_{ z\in [\frac{b+x}{2},b)}
\frac{\gam_2(z-a)(b-z)}{\hat{\kappa} \kappa (b-a)|z-x|}
\omega_{y,z}\Big\}.
\end{align*}
Now, observe that
\begin{align}
\label{Gut1}
\frac{(x-a)(b-x)}{b-a}\sum_{ z\in (a,\frac{x+a}{2}]}
\frac{(z-a)(b-z)}{ (b-a)|z-x|} \omega_{y,z}
&\leq \sum_{ z\in (a,\frac{x+a}{2}]}(b-z)\omega_{y,z}\nonumber\\
&\leq \sum_{z<b}(b-z)\omega_{y,z},
\end{align}

\begin{align}
\label{Gut2}
\frac{(x-a)(b-x)}{b-a}\sum_{ z\in (\frac{x+a}{2},
\frac{b+x}{2})} \omega_{y,z}
&\leq (b-x)\sum_{ z\in (\frac{x+a}{2},\frac{b+x}{2})} 
 \omega_{y,z}\nonumber\\
&\leq 2\sum_{ z\in (\frac{x+a}{2},\frac{b+x}{2})} (b-z) 
\omega_{y,z}\nonumber\\
&\leq 2 \sum_{ z<b} (b-z) \omega_{y,z}
\end{align}
and
\begin{align}
\label{Gut3}
\frac{(x-a)(b-x)}{b-a}\sum_{ z\in [\frac{x+b}{2},b)}
\frac{(z-a)(b-z)}{ (b-a)|z-x|} \omega_{y,z}
&\leq 2\sum_{ z\in [\frac{x+b}{2},b)}(b-z)\omega_{y,z}\nonumber\\
&\leq 2\sum_{z<b}(b-z)\omega_{y,z}.
\end{align}
Putting~(\ref{Gut1}), (\ref{Gut2}) and~(\ref{Gut3}) together leads
to 
\begin{equation}
H_1(M)\leq 2\Big( \frac{2\gam_2}{\hat{\kappa} \kappa}+1\Big) 
\sum_{y\in J_1}\sum_{z<b}(b-z)\omega_{y,z}.
\end{equation}
Observe that this last upper bound on~$H_1(M)$ does not depend
on~$x$ anymore. Now, by Condition~K, for any $\eta>0$, we can
take $M_1>0$ sufficiently large such that $\IP$-a.s.\ for all $u\in
\Z$ we have
\[
\sum_{v>u+M_1}\sum_{w< u}\omega_{v,w}(u-w)
 <\frac{ \kappa\eta}{4} \Big(
\frac{2\gam_2}{\hat{\kappa} \kappa}+1\Big)^{-1}.
\]
For this~$M_1$, we have $H_1(M_1)< \kappa\eta/2$.
By symmetry, we also have that $H_2(M_1)< \kappa\eta/2$.
Combining these two last results with (\ref{Simplif}) and (\ref{Ertac}) and the case of intervals of length 2 treated at the beginning of the proof, we obtain 
that for every $\eta>0$ there exists $M$ 
such that $\IP$-a.s., for any interval $[a,b]$,
\begin{equation*}
\min_{x\in (a,b)}\Po^x[X_{\tau_{E}}\in I_M]\geq 1-\eta.
\end{equation*}
This concludes the proof of Proposition~\ref{SNL}. 
\qed

\section{Proof of Theorem~\ref{Theouni}}
\label{s_proof_t_uni}
In this section we prove the UCLT. Let $\C^u_b(C(\R_+),\R)$ be the space of bounded
uniformly continuous functionals from $C(\R_+)$ into $\R$. First, let us prove the apparently weaker statement:
\begin{prop}
\label{UCLTprop1}
For all $F\in \C^u_b(C(\R_+),\R)$, we have $\IP$-a.s., for every $H>0$, 
\begin{equation*}
\lim_{n\to \infty}\sup_{x\in [-H\sqrt{n},H\sqrt{n}]}\Big| 
{\mathtt E}_{\theta_x \omega}[F(Z^{n})]-E[F(W)]\Big|=0.
\end{equation*}
\end{prop}
The difficult part of the proof of Theorem~\ref{Theouni} is to show 
Proposition~\ref{UCLTprop1}. To prove this proposition, we will introduce the
notion of ``good site" in $\Z$. The set of good sites is by definition
the set of sites in~$\Z$ from which we can guarantee that the random walk
converges uniformly to Brownian motion. Due to the ergodicity of the random
environment, we will then prove that starting a random walk from any site
in $[-H\sqrt{n},H\sqrt{n}]$, with high probability, it will meet a close good
site quickly enough to derive a uniform CLT.  This part will be done in
two steps, introducing the intermediate concept of ``nice site".  
More precisely, the sequence of steps we will follow in  this section to prove
Proposition~\ref{UCLTprop1} is the following:
\begin{itemize}
\item In Definition~\ref{goodsite}, we formally define the notion of ``good
sites''.
\item In Definition~\ref{nicesite}, we introduce the notion of ``nice
sites".  Heuristically, $x$ is a nice site if for some $\delta>0$ and
$h>0$, the range of the random walk starting from $x$ until time $hn$ is
greater than $\delta h^{1/2} \sqrt{n}$ with high probability, so that the random
walk cannot stay ``too close'' to its starting location (it holds that good
sites are nice). 
\item Right after Definition~\ref{nicesite}, 
we show that any interval $I\in [-2H\sqrt{n}, 2H\sqrt{n}]$ of 
length $n^{\nu}$ with $\nu\in (1/(2+\beta),\frac{1}{2})$ 
(here $\beta$ is from Condition~K) must contain at least one ``nice site".
\item In Lemma~\ref{reacnicesite} we show that,
 starting from a site $x\in [-H\sqrt{n}, H\sqrt{n}]$, with high probability the random walk
meets a nice site at a distance at most $n^{\mu}$ in time at most  $n^{2\mu}$
with $\mu \in (\nu,\frac{1}{2})$.
\item In Lemma~\ref{reacgoodsite} we show that,
starting from a nice site $x\in [-(3/2)H\sqrt{n}, (3/2)H\sqrt{n}]$, the random walk meets
with high probability a good site at a distance less than $h\sqrt{n}$ before
time $hn$. 
\item We combine Lemmas~\ref{reacnicesite} 
and~\ref{reacgoodsite} to obtain that, starting from any $x\in [-H\sqrt{n},H\sqrt{n}]$ the
random walk meets a good site at a distance less than $h\sqrt{n}$ before time~$hn$. 
This is the statement of Lemma~\ref{reacgoodsitefrombad}. 
\item Proposition~\ref{UCLTprop1} then follows from Lemma~\ref{reacgoodsitefrombad}, since we know essentially 
that the random walk will quickly reach a nearby good site, and
 from this good site the convergence properties are ``good'' by definition. 
\end{itemize}

 From Proposition~\ref{UCLTprop1}, 
we obtain Theorem~\ref{Theouni} in the following way. 
In Proposition~\ref{portmantpartial}, we first show that Proposition~\ref{UCLTprop1} implies a corresponding statement in which we substitute uniformly continuous functionals $F$ by open sets of $C(\R_+)$. Then, in Proposition~\ref{portlast}, we use the separability of the space $C(\R_+)$ to show that we can interchange the terms ``for any open set $G$" and ``$\IP$-a.s." in (ii) of Proposition~\ref{portmantpartial}. 
Then, we use standard arguments as in the proof of 
the Portmanteau theorem of~\cite{Billing} to conclude the proof of Theorem~\ref{Theouni}.

\medskip
 
Now, fix $F\in \C^u_b(C(\R_+),\R)$. 
Our first goal is to prove Proposition~\ref{UCLTprop1}, that is, $\IP$-a.s., for every
$\tilde{\eps}, H>0$, 
\begin{equation}
\label{UCLTrefor}
\sup_{x\in [-H\sqrt{n},H\sqrt{n}]}\Big| 
{\mathtt E}_{\theta_x \omega}[F(Z^{n})]-E[F(W)]\Big|\leq
\tilde{\eps}
\end{equation}
for all large enough~$n$. To start, we need to 
write some definitions and prove some intermediate results. 
From now on, we suppose that $\sig=1$ 
(otherwise replace~$X$ by $\sig^{-1}X$).

Denote
\begin{align*}
 R_n^+(m) &= \max_{s\leq m}(X_{n+s}-X_n),\\
 R_n^-(m) &= \min_{s\leq m}(X_{n+s}-X_n),\\
 R_n(m) &=R_n^+(m) - R_n^-(m),
\end{align*}
and 
\begin{align*}
{ \mathfrak R}_t^+(u) &= \max_{s\leq u}(W_{t+s}-W_t),\\
 { \mathfrak R}_t^-(u) &= \min_{s\leq u}(W_{t+s}-W_t),\\
 { \mathfrak R}_t(u) &={ \mathfrak R}_t^+(u) 
  - { \mathfrak R}_t^-(u).
\end{align*}
Let~$\dd$ be the distance on the space $C_{\R_+}$ 
defined by
\[
\dd(x,y)=\sum_{n=1}^{\infty}2^{-n+1}
\min\Big\{1,\sup_{s\in [0,n]}|x(s)-y(s)|\Big\}.
\]
Now, for any given $\eps>0$, we define
\begin{align}
\label{delchos}
\delta_\eps&:=\max\Big\{\delta_1\in (0,1]:P[{\mathfrak R}_0(1/2)< \delta_1]+P[{\mathfrak R}_{1/2}(1/2)
 < \delta_1]\nonumber\\
 &\phantom{******}+P [{\mathfrak R}_{1}(1/2)< \delta_1]+P [ {\mathfrak
R}^{+}_0(1)< \delta_1]+P [ {\mathfrak
R}^{-}_0(1)< \delta_1]\leq \frac{\eps}{2}\Big\}
\end{align}
and
\begin{equation}
\label{UUFMW}
h_\eps:=\max\Big\{h_1\in (0,1]:P\Big[\sup_{s\leq h_1}|W(s)|>\eps\Big]
+P\Big[\sup_{s\leq h_1}\dd(\theta_sW,W)>\eps\Big]\leq
\frac{\eps}{2}\Big\}.
\end{equation}
Observe that $\delta_\eps$ and $h_\eps$ are positive for all $\eps>0$ and decrease to 0 as $\eps\to 0$. For~(\ref{UUFMW}), the positivity of  $h_\eps$ for $\eps>0$ follows from the properties of the modulus of continuity of Brownian motion (see e.g.\ Theorem~1.12 of~\cite{PerMot}).
\begin{df}
\label{goodsite}
For a given realization $\omega$ of the environment and $N\in \N$, we say that $x\in\Z$ is $(\eps,N)$-good,
 if 
\begin{itemize}
\item[(i)]
$\min \Big \{n\geq 1: \big| 
{\mathtt E}_{\omega}[F(Z^{m})]-E[F(W)]\big| \leq \eps,
\phantom{*}\mbox{for all $m\geq n$}\Big \}\leq N$;
\item[(ii)]
$ \Po^x\Big[R_k(h_\eps m)\geq \delta_\eps h_\eps^{1/2}\sqrt{m}\text{ for all }
k\leq h_\eps m, R_0^\pm(h_\eps m)\geq\delta_\eps h_\eps^{1/2}\sqrt{m}\Big] \geq 1-\eps,$
for all $m\geq N$;
\item[(iii)]
$ \mathtt{P}_{\theta_x \omega}\Big[\sup_{s\leq h_\eps}|Z^{m}(s)|
\leq \eps,\sup_{s\leq
h_\eps}\dd(\theta_sZ^{m},Z^{m})\leq \eps\Big] 
\geq 1-\eps$, for
all $m\geq N$.
\end{itemize}
\end{df}

For any given $\eps>0$, it follows from Theorem~\ref{t_q_invar_princ}, (\ref{delchos}) and (\ref{UUFMW}) that for any $\eps'>0$ there exists~$N$ such that
\[
 \IP[0\text{ is $(\eps,N)$-good}] 
        > 1-\eps'.
\]
 Then, by the Ergodic Theorem, $\IP$-a.s., for all $n$ large enough, it holds that 
\begin{equation}
\label{number_of_not_good}
 \big|\{x\in [-2H\sqrt{n},2H\sqrt{n}] : x\text{ is not $(\eps,N)$-good}\}\big| < 5\eps'H\sqrt{n}.
\end{equation}

Next, we need the following
\begin{df}
\label{nicesite}
We say that a site~$x$ is $(\eps, n)$-nice, if 
\[
 \Po^x\Big[R_0(h_\eps n)\geq\delta_\eps h_\eps^{1/2}\sqrt{n}\Big]\geq 1-3\eps.
\]
\end{df}
In particular, note that, if for some~$N\leq n$ a site~$x$ 
is $(\eps,N)$-good,
then it is $(\eps, n)$-nice. 

Now, fix some $\nu \in(\frac{1}{2+\beta},\frac{1}{2})$ (so that
$\nu(2+\beta)-1>0$), where $\beta$ is from Condition~K.
Observe that by Condition~K there exists $\gam_1>0$ such that 
for any starting point~$x\in\Z$
\begin{equation}
\label{eq_long_jump}
 \Po^x[|X_{k+1}-X_k|<n^{\nu} \text{ for all }k\leq h_\eps n]
 \geq 1 - \gam_1n^{-(\nu(2+\beta)-1)}.
\end{equation}
We argue by contradiction that, $\IP$-a.s., there exists $n_1=n_1(\omega, \eps)$ such that any interval of length at least $n^{\nu}$ contains a $(\eps,n)$-nice site for $n>n_1$. For this, choose $\eps'>0$ such that $5\eps'H<\delta_\eps h_\eps^{1/2}$ and let $n$ large enough such that
$n^\nu<5\eps'H\sqrt{n}$ and (\ref{number_of_not_good}) hold. Let $I\subset [-2H\sqrt{n},2H\sqrt{n}]$ be an
interval of length~$n^\nu$ such that it does not contain
any $(\eps, n)$-nice site. Observe that,
by~\eqref{number_of_not_good}, there exists a
$(\eps,N)$-good site~$x_0$ such that 
$|x_0-y|<\delta_\eps h_\eps^{1/2}\sqrt{n}$ for all $y\in I$. Note that,
by~\eqref{eq_long_jump} and (ii) of Definition \ref{goodsite},
\[
 \Po^{x_0}[\text{there exists }k\leq h_\eps n\text{ such that }X_k\in
I] \geq 1-\eps-\gam_1n^{-(\nu(2+\beta)-1)}
\]
(the particle crosses~$I$ without jumping over it entirely), hence
\begin{align*}
\Po^{x_0}[\text{there exists }k\leq h_\eps n\text{ such that }
  R_k(h_\eps n)<\delta_\eps h_\eps^{1/2}\sqrt{n}]
 & \geq 3\eps (1-\eps-\gam_1n^{-(\nu(2+\beta)-1)})\\
 & \geq 2\eps
\end{align*}
if~$n$ is large enough. But this contradicts the fact that $x_0$ is 
$(\eps,N)$-good.
So, we see that, $\IP$-a.s., any interval $I\subset [-2H\sqrt{n},2H\sqrt{n}]$ of
length~$n^\nu$ should contain at least one 
$(\eps,n)$-nice site for $n$ large enough. 

Let $\mu \in (\nu,\frac{1}{2})$. In the next Lemma, we show that starting from a 
site $x\in
[-H\sqrt{n},H\sqrt{n}]$, with high probability the random walk will meet a
$(\eps,n)$-nice site at a distance at most~$n^{\mu}$ in 
time at most~$n^{2 \mu}$. 

For $x\in \Z$ and $n\in\N$, let us denote by
$\mathcal{P}^l_n(x)$ the largest $y\leq x$ such that~$y$ is a
$(\eps,n)$-nice site and by $\mathcal{P}^r_n(x)$ the smallest
$y\geq x$ such that $y$ is a $(\eps,n)$-nice site. Furthermore, we
denote by $\mathcal{N}_{\eps,n}$ the set of $(\eps, n)$-nice
sites in $\Z$.
\begin{lm}
\label{reacnicesite}
 For any $\eps_1>0$ and $\eps>0$, we have $\IP$-a.s., for all sufficiently large $n$, for all $x\in [-H\sqrt{n},H\sqrt{n}]$,
\begin{equation}
\label{TEREN}
\Po^x\Big[\tau_{\mathcal{N}_{\eps,n}}\leq n^{2 \mu},
 \max_{j\leq \tau_{\mathcal{N}_{\eps,n}}}|X_{j}-X_0|\leq
n^{\mu}\Big]\geq 1-\eps_1.
\end{equation}
\end{lm}
\textit{Proof.} First, suppose that $x$ is not $(\eps,n)$-nice, otherwise the proof of (\ref{TEREN}) is trivial.
For some integer $M_1>0$, define the intervals 
\[
I_{M_1}(x):=[\mathcal{P}^l_n(x)-M_1, \mathcal{P}^r_n(x)+M_1].
\]
Let us also define the following increasing sequence of stopping times: 
$\xi_0:=0$ and for $i\geq 1$,
\[
\xi_i:=\inf\big\{k>\xi_{i-1}: X_k\notin 
(\mathcal{P}^l_n(X_{\xi_{i-1}}),\mathcal{P}^r_n(X_{\xi_{i-1}}))\big\}
+M_1.
\]
Then, we define the events 
\begin{align*}
 A_i&:=\big\{\mbox{there exists $k \in (\xi_i, \xi_{i+1}]$ such that 
$X_k\in \{\mathcal{P}^l_n(X_{\xi_{i}}),
\mathcal{P}^r_n(X_{\xi_{i}})$}\}\big\},\\
B_i &:=\Big\{X_{\xi_{i+1}-M_1}\in I_{M_1}(X_{\xi_i}), 
\max_{(k,l)\in [\xi_{i+1}-M_1, \xi_{i+1}]^2}|X_k-X_l|\leq
M_1n^{\nu}\Big\}
\end{align*}
for $i\geq 0$.
Now, fix temporarily $\tilde{\eps}\in (0,\frac{1}{2})$. By Proposition~\ref{SNL}, we can choose~$M_1$ large enough
in such a way that 
\begin{equation}
\label{CGH}
\min_{y\in (\mathcal{P}^l_n(x),\mathcal{P}^r_n(x))}
\Po^y[X_{\xi_1-M_1}\in I_{M_1}(x)]\geq 1-\tilde{\eps}.
\end{equation}
Note that by Condition~E and Proposition~\ref{SNL} and (\ref{CGH}) we have 
$\Po[A_0]\geq \frac{1}{2}\kappa^{M_1}$.
By the Markov property, we obtain for some integer $L>0$,
\begin{equation}
\label{DROT1}
\Po^x\Big[\bigcap_{i=0}^{L-1}A_i^c\Big]
=\Po^x[A_0^c]\dots\Po^x[A_{L-1}^c\mid A_0^c\dots A_{L-2}^c]\leq
\Big(1-\frac{1}{2}\kappa^{M_1}\Big)^L.
\end{equation}
Since the event $\{\max_{(k,l)\in [\xi_{1}-M_1, \xi_{1}]^2}|X_k-X_l|>
M_1n^{\nu}\}$ implies that there is at least one jump of size $n^{\nu}$ during the time interval $[\xi_{1}-M_1, \xi_{1}]$, using Condition K, (\ref{Ellip2}), (\ref{CGH}) and the Markov property, we have that $\Po^x[B_0^c]\leq \tilde{\eps}+\gam_1n^{-\nu(2+\beta)}$ for some positive constant $\gam_1$.
We obtain by the Markov property,
\begin{equation}
\label{DROT2}
\Po^x\Big[\bigcup_{i=0}^{L-1}B_i^c\Big]
\leq L(\tilde{\eps}+\gam_1n^{-\nu(2+\beta)}).
\end{equation}
Observe that each event $\{\max_{j\in [\xi_i-M_1, \xi_i]}|X_{j}-X_{\xi_i-M_1-1}|>
(M_1+1)n^{\nu}+M_1\}$, $ i\geq 1$, implies either $\{|X_{\xi_i-M_1}-X_{\xi_i-M_1-1}|>n^{\nu}+M_1\}$ (which implies that the first jump after time $\xi_i-M_1-1$ is out of the interval $I_{M_1}(\xi_i-M_1-1)$) or $\{\max_{(k,l)\in [\xi_{i}-M_1, \xi_{i}]^2}|X_k-X_l|>
M_1n^{\nu}\}$. Then, combining~(\ref{DROT1}) and~(\ref{DROT2}), we have
\begin{align}
\label{DROT3}
\lefteqn{
\Po^x\Big[\tau_{\mathcal{N}_{\eps,n}} \in [0, \xi_L], 
\max_{j\leq \tau_{\mathcal{N}_{\eps,n}}}|X_{j}-X_0|\leq
L((M_1+1)n^{\nu}+M_1)\Big]
}\nonumber\\
&\phantom{*****************}
\geq
1-\Big(1-\frac{1}{2}\kappa^{M_1}\Big)^L-L\tilde{\eps}-L\gam_1n^{
-\nu(2+\beta)}.
\end{align}
Now, let $\mu'>0$ and denote 
$G_i:=\{\xi_i-\xi_{i-1}\leq n^{2\nu+\mu'}+M_1 \}$
for $1\leq i\leq L$.
We have
\begin{align}
\label{Poit}
\Po^x[\xi_L\leq L(n^{2\nu+\mu'}+M_1)]
&\geq \Po^x[G_1,G_2,\dots,G_L]\nonumber\\
&= \Po^x[G_1]\Po^x[ G_2\mid G_1]\dots 
\Po^x[G_L\mid G_1\dots G_{L-1}].
\end{align}
By Proposition~\ref{specgap} and the fact that any interval of
length~$n^{\nu}$ in $[-2H\sqrt{n}, 2H\sqrt{n}]$ should contain at least one
$(\eps,n)$-nice site, we obtain 
\begin{equation*}
\Po^x[G_1]\geq 1-\exp\Big(-\frac{n^{\mu'}}{K_1}\Big)
\end{equation*}
for sufficiently large $n$.
By~(\ref{Poit}) and the Markov property we have
\begin{align}
\label{DROT4}
\Po^x[\xi_L\leq L(n^{2\nu+\mu'}+M_1)]
\geq \Big[1-
 \exp\Big(-\frac{n^{\mu'}}{ K_1}\Big)\Big]^L
\geq 1-L\exp\Big(-\frac{n^{\mu'}}{K_1}\Big).
\end{align}
Now, choose~$L$ sufficiently large so that 
\begin{equation}
\label{DROT5}
\Big(1-\frac{1}{2}\kappa^{M_1}\Big)^L\leq 
\frac{\eps_1}{3}
\end{equation}
and $\tilde{\eps}$ sufficiently small such that
\begin{equation}
\label{DROT6}
L\tilde{\eps}\leq \frac{\eps_1}{3}.
\end{equation}
Then, combining~(\ref{DROT3}), (\ref{DROT4}), 
 (\ref{DROT5}) and~(\ref{DROT6}), we obtain 
\begin{align*}
\lefteqn{\Po^x\Big[\tau_{\mathcal{N}_{\eps, n}}
\leq L(n^{2\nu+\mu'}+M_1),\max_{j\leq
\tau_{\mathcal{N}_{\eps,n}}}|X_{j}-X_0|\leq
L((M_1+1)n^{\nu}+M_1)\Big]}\phantom{********************}\nonumber\\
&\geq 1-\frac{2\eps_1}{3}-L\gam_1n^{-\nu(2+\beta)}-L
\exp\Big(-\frac{n^{\mu'}}{ K_1}\Big).
\end{align*}
Taking $\mu'<2(\mu-\nu)$ we obtain for all sufficiently 
large~$n$
\[
\Po^x\Big[\tau_{\mathcal{N}_{\eps,n}}\leq n^{2 \mu},
\max_{j\leq \tau_{\mathcal{N}_{\eps,n}}}|X_{j}-X_0|\leq
n^{\mu}\Big]\geq 1-\eps_1.
\]
This concludes the proof of Lemma~\ref{reacnicesite}. 
\qed

\medskip

 We now show that we can find $\eps>0$ small enough such that starting from a $(\eps,n)$-nice 
site $x\in [-\frac{3}{2}H\sqrt{n},\frac{3}{2}H\sqrt{n}]$, with high probability, the random walk will meet a
$(\eps,n)$-good site 
 at a distance at most $h_\eps^{1/2}\sqrt{n}$ before time~$h_{\eps}n$. We denote by $\mathcal{G}_{\eps,N}$ the set of $(\eps,N)$-good sites in $\Z$.
\begin{lm}
\label{reacgoodsite}
 For any $\eps_1>0$ and $\eps\in (0, \frac{\eps_1}{6}]$, we have that $\IP$-a.s., for all sufficiently large~$n$, 
for all $x\in [-\frac{3}{2}H\sqrt{n},\frac{3}{2}H\sqrt{n}]\cap \mathcal{N}_{\eps,n}$,
\[
\Po^x\Big[\tau_{\mathcal{G}_{\eps,n}}\leq h_\eps n,
 \max_{j\leq \tau_{\mathcal{G}_{\eps,n}}}|X_{j}-X_0|\leq
h_\eps^{1/2}\sqrt{n}\Big]\geq 1-\eps_1.
\]
\end{lm}

\noindent\textit{Proof.}
Fix some integer $M>1$ and consider the following partition 
of~$\Z$ into intervals of size~$M$:
\[
J_j=[jM,(j+1)M),\phantom{**}j\in \Z.
\]
We say that an interval $J_j$ is $(\eps,N)$-good if all the points inside~$J_j$ are
$(\eps,N)$-good (otherwise we call the interval ``bad").
Fix $\eps\leq \frac{\eps_1}{6}$. Then,we have that, for any $\eps'>0$,
there exists~$N$ such that
\[
 \IP[J_0\text{ is $(\eps,N)$-good}] 
      > 1-\eps'.
\]
Then, by the Ergodic
Theorem, $\IP$-a.s., for all $n$ large enough it holds that 
\begin{equation}
\label{number_of_not_good1}
 \big|\{J_j\text{ such that } j\in [-H\sqrt{n},H\sqrt{n}] 
\text{ and $J_j$ is not
$(\eps,N)$-good}\}\big| <
3\eps'H\sqrt{n}.
\end{equation}
In particular, from this last inequality, we deduce that the 
length of the largest subinterval of $[-2H\sqrt{n}, 2H\sqrt{n}]$ that is a union of bad intervals is
smaller than $3\eps' HM\sqrt{n}$. Let $x\in[-\frac{3}{2}H\sqrt{n},\frac{3}{2}H\sqrt{n}]$ be a 
$(\eps,n)$-nice
site that belongs to an interval 
$I=[a,b]\subset [-2H\sqrt{n},2H\sqrt{n}]$ that is a maximal union of bad $J_j$'s (so that the adjacent $J_j$'s to $I$ are necessarily good).
Then, choose~$\eps'$ such that $ 3\eps'HM<h^{1/2}\delta_\eps$. Thus, 
in a time of order~$h_\eps n$ a random walk starting at~$x$ will leave the interval $I$ with high probability. When this happens, to guarantee that the random walk
will hit a $(\eps,N)$-good site with high
probability, we can choose a large enough~$M$ in such a way that 
\begin{equation}
\label{RTT1}
\Po^x[X_{\tau_{I^c}}\in I_M]\geq 1 -\frac{\eps_1}{2},
\end{equation}
with $I_M=[a-M,a]\cup[b,b+M]$.
By definition of a $(\eps,n)$-nice site, 
since $\eps\leq \frac{\eps_1}{6}$ we have 
\begin{equation}
\label{RTT2}
\Po^x[\tau_{I^c}\leq h_\eps n]\geq \Po^x[R_0(h_\eps n)\geq \delta_\eps h_\eps^{1/2} \sqrt{n}]
\geq 1-3\eps\geq 1-\frac{\eps_1}{2}.
\end{equation}
Thus, combining~(\ref{RTT1}) and~(\ref{RTT2}), and using the fact
that $\delta_\eps\in (0,1]$, we obtain $\IP$-a.s.,
\begin{equation*}
\Po^x\Big[\tau_{\mathcal{G}_{\eps,n}}\leq h_\eps n,
\max_{j\leq \tau_{\mathcal{G}_{\eps,n}}}|X_{j}-X_0|\leq
h_\eps^{1/2}\sqrt{n}\Big]\geq 1-\eps_1.
\end{equation*}
 for all~$n$ large enough and $x\in [-\frac{3}{2}H\sqrt{n},\frac{3}{2}H\sqrt{n}]\cap \mathcal{N}_{\eps,n}$.
This concludes the proof of Lemma~\ref{reacgoodsite}. 
\qed

Then, combining Lemmas~\ref{reacnicesite} 
and~\ref{reacgoodsite}, we can deduce (considering for example rational values for $\eps_1$ and $\eps$):
\begin{lm}
\label{reacgoodsitefrombad}
The following statement holds $\IP$-a.s.: for any $\eps_1>0$, we can choose~$\eps>0$ arbitrary small in such a way that for all sufficiently large~$n$ and for all $x\in [-H\sqrt{n},H\sqrt{n}]$,
\[
\Po^x\Big[\tau_{\mathcal{G}_{\eps,n}}\leq h_\eps n,
\max_{j\leq \tau_{\mathcal{G}_{\eps,n}}}|X_{j}-X_0|\leq
h_\eps^{1/2}\sqrt{n}\Big]\geq 1-\eps_1.
\]
\end{lm}

We are now ready to prove Proposition~\ref{UCLTprop1}.
\medskip

\noindent
\textit{Proof of Proposition~\ref{UCLTprop1}.}
Let us prove (\ref{UCLTrefor}). Let $x\in [-H\sqrt{n},H\sqrt{n}]$. 
In this
last part, for the sake of brevity, we denote by $\G$ the set of
$(\eps,n)$-good sites.
 Let us denote by
\[
R:=\Big| {\mathtt E}_{\theta_x\omega}[F(Z^{n})]-E[F(W)]\Big|
\]
the quantity we have to bound.
Let $\eps \leq \frac{\tilde{\eps}}{2}$, we have
by definition of a $(\eps,n)$-good site,
\begin{align}
\label{UFin}
R
&\leq \Big|{\mathtt E}_{\theta_x\omega}\Big(\!F(Z^{n})
- {\mathtt E}_{\theta_{X_{\tau_{\G}}} \omega}[F(Z^{n})]\Big)\Big|
+ \Big|{\mathtt E}_{\theta_x\omega}
\Big({\mathtt E}_{\theta_{X_{\tau_{\G}}}
\omega}[F(Z^{n})]-E[F(W)]\Big)\Big|\nonumber\\
&\leq \Big|{\mathtt E}_{\theta_x\omega}\Big(\!F(Z^{n})
- {\mathtt E}_{\theta_{X_{\tau_{\G}}} \omega}[F(Z^{n})]\Big)\Big| + \frac{\tilde{\eps}}{2}.
\end{align}
Denote $X':=X-x$ and observe that, by the Markov property
\begin{align}
\label{ZION-1}
\Big|{\mathtt E}_{\theta_x\omega}\Big(\!F(Z^{n})
- {\mathtt E}_{\theta_{X_{\tau_{\G}}} \omega}[F(Z^{n})]\Big)\Big|
&= \Big|{\mathtt E}_{\theta_x\omega}\Big(F(Z^{n})-{\mathtt
E}_{\theta_{X'_{\tau_{\G}}}( \theta_x
\omega)}[F(Z^{n})]\Big)\Big|\nonumber\\
&=\Big|{\mathtt E}_{\theta_x\omega}[F\circ Z^{n}
 -F\circ \theta_{n^{-1}\tau_\G}(
Z^{n}-n^{-\frac{1}{2}}X'_{\tau_\G})]\Big|\nonumber\\
&\leq {\mathtt E}_{\theta_x \omega}
\Big|F\circ Z^{n}-F\circ \theta_{n^{-1}\tau_\G}(
Z^{n}-n^{-\frac{1}{2}}X'_{\tau_\G})\Big|.
\end{align}
We are going to show that for all sufficiently large~$n$ we 
have uniformly in $x\in [-H\sqrt{n},H\sqrt{n}]$
\[
{\mathtt E}_{\theta_x\omega}\Big|F\circ Z^{n}
 -F\circ \theta_{n^{-1}\tau_\G}(
Z^{n}-n^{-\frac{1}{2}}X'_{\tau_\G})\Big|\leq \frac{\tilde{\eps}}{2}
\]
for $\eps>0$ small enough.
Let $M^{n}:=Z^{n}-n^{-\frac{1}{2}}X'_{\tau_\G}$.
Since~$F$ is uniformly continuous, we can choose $\eta>0$ in such
a way that if $\dd(x,y)\leq\eta$ then $|F(x)-F(y)|\leq
\frac{\tilde{\eps}}{4}$. Then, we have
\begin{align}
\label{ZION0}
{\mathtt E}_{\theta_x\omega} \Big|F\circ
Z^{n}-F\circ\theta_{n^{-1}\tau_{\mathcal{G}}}M^{n} \Big|
&={\mathtt E}_{\theta_x\omega}
\Big[ \Big|F\circ
Z^{n}-F\circ\theta_{n^{-1}\tau_{\mathcal{G}}}M^{n} \Big|\1{\dd(Z^{n},\theta_{n^{-1}\tau_{\mathcal{G}}}M^{n})\leq \eta}\Big]\nonumber\\
&\phantom{**}
 +{\mathtt E}_{\theta_x\omega}\Big[\Big|F\circ
Z^{n}-F\circ\theta_{n^{-1}\tau_{\mathcal{G}}}M^{n} \Big|\1{\dd(Z^{n},\theta_{n^{-1}\tau_{\mathcal{G}}}M^{n})> \eta}\Big]\nonumber\\
&\leq \frac{\tilde{\eps}}{4}+2\|F\|_{\infty}
{\mathtt
P}_{\theta_x\omega}\Big[\dd(Z^{n},
\theta_{n^{-1}\tau_{\mathcal{G}}}M^{n})>\eta\Big].
\end{align}
Since $h_\eps \leq 1$, we have
\begin{align}
\label{ZION}
{\mathtt P}_{\theta_x\omega}
\Big[\dd(Z^{n},\theta_{n^{-1}\tau_{\mathcal{G}}}M^{n})>
\eta\Big]
&\leq {\mathtt P}_{\theta_x\omega}\Big[\dd(Z^{n},
\theta_{n^{-1}\tau_{\mathcal{G}}}M^{n})>\eta,
\tau_{\mathcal{G}}\leq h_\eps n\Big]+{\mathtt
P}_{\theta_x\omega}[\tau_{\mathcal{G}}> h_\eps n]\nonumber\\
&\leq {\mathtt P}_{\theta_x\omega}
\Big[\sup_{t\in [0, n^{-1}\tau_\G]}|
Z^{n}-\theta_{n^{-1}\tau_{\mathcal{G}}}M^{n}|>\frac{\eta}{2},
\tau_{\mathcal{G}}\leq h_\eps n\Big]\nonumber\\
&\phantom{**}+{\mathtt P}_{\theta_x\omega}
\Big[\dd(\theta_{n^{-1}\tau_{\mathcal{G}}}
Z^{n},\theta^2_{n^{-1}\tau_{\mathcal{G}}}M^{n})>\frac{\eta}{2},
\tau_{\mathcal{G}}\leq h_\eps n\Big]+{\mathtt
P}_{\theta_x\omega}[\tau_{\mathcal{G}}> h_\eps n].\nonumber\\
\end{align}
Let $\mathcal{F}_{\tau_\G}$ be the $\sigma$-field generated by~$X$
until time $\tau_\G$. We decompose the first term in the
right-hand side of~(\ref{ZION}) in the following way
\begin{align}
\label{ZION1}
\lefteqn{
{\mathtt P}_{\theta_x\omega}
\Big[\sup_{t\in [0, n^{-1}\tau_\G]}|
Z^{n}-\theta_{n^{-1}\tau_{\mathcal{G}}}M^{n}|>\frac{\eta}{2},
\tau_{\mathcal{G}}\leq
h_\eps n\Big]
}\phantom{****}\nonumber\\
&\leq {\mathtt P}_{\theta_x\omega}
\Big[\sup_{t\in [0, n^{-1}\tau_\G]}|
Z^{n}|>\frac{\eta}{4}\Big]+{\mathtt
P}_{\theta_x\omega}\Big[\sup_{t\in [0,
h_\eps]}|\theta_{n^{-1}\tau_{\mathcal{G}}}M^{n}|>
\frac{\eta}{4}\Big]\nonumber\\
&= {\mathtt P}_{\theta_x\omega}
\Big[\sup_{t\in [0, n^{-1}\tau_\G]}|
Z^{n}|>\frac{\eta}{4}\Big]
+{\mathtt
E}_{\theta_x\omega}\Big(\mathtt{P}_{\theta_x\omega}\Big[\sup_{t\in
[0, h_\eps]}|\theta_{n^{-1}\tau_{\mathcal{G}}}M^{n}|>
\frac{\eta}{4} \;\Big|\; \mathcal{F}_{\tau_\G}
\Big]\Big)\nonumber\\
&= {\mathtt P}_{\theta_x\omega}
\Big[\sup_{t\in [0, n^{-1}\tau_\G]}|
Z^{n}|>\frac{\eta}{4}\Big]+{\mathtt
E}_{\theta_x\omega}\Big(\mathtt{P}_{\theta_{X_{\tau_\G}}\omega}\Big[
\sup_{t\in [0, h_\eps]}|Z^{n}|>\frac{\eta}{4}\Big]\Big).
\end{align}
We now deal with the second term of the right-hand side
of~(\ref{ZION})
\begin{align}
\label{ZION2}
\lefteqn{{\mathtt P}_{\theta_x\omega}
\Big[\dd(\theta_{n^{-1}\tau_{\mathcal{G}}}
Z^{n},\theta^2_{n^{-1}\tau_{\mathcal{G}}}M^{n})>
\frac{\eta}{2},\tau_{\mathcal{G}}\leq
h_\eps n\Big]}\phantom{***}\nonumber\\
&\leq {\mathtt P}_{\theta_x\omega}\Big[|X'_{\tau_\G}|>\frac{\eta}{4}\sqrt{n}\Big]+
{\mathtt
P}_{\theta_x\omega}\Big[\dd(\theta_{n^{-1}\tau_{\mathcal{G}}}
M^{n},\theta^2_{n^{-1}\tau_{\mathcal{G}}}M^{n})>
\frac{\eta}{4},\tau_{\mathcal{G}}\leq h_\eps n\Big]\nonumber\\
&\leq  {\mathtt P}_{\theta_x\omega}\Big[\sup_{t\in [0,
n^{-1}\tau_\G]}| Z^{n}|>\frac{\eta}{4}\Big]
+{\mathtt E}_{\theta_x\omega}
\Big(\1{\tau_{\mathcal{G}}\leq h_\eps n}{\mathtt
P}_{\theta_x\omega}\Big[\dd(\theta_{n^{-1}\tau_{\mathcal{G}}}
M^{n},\theta^2_{n^{-1}\tau_{\mathcal{G}}}M^{n})>
\frac{\eta}{4} \;\Big|\; \mathcal{F}_{\tau_\G}\Big]\Big)\nonumber\\
&= {\mathtt P}_{\theta_x\omega}
\Big[\sup_{t\in [0, n^{-1}\tau_\G]}|
Z^{n}|>\frac{\eta}{4}\Big]
+{\mathtt E}_{\theta_x\omega}
\Big(\1{\tau_{\mathcal{G}}\leq
h_\eps n}\mathtt{P}_{\theta_{X_{\tau_\G}}\omega}
\Big[\dd(Z^{n},\theta_{n^{-1}\tau_{\mathcal{G}}}Z^{n})
>\frac{\eta}{4}\Big]\Big).
\end{align}
Combining~(\ref{ZION}), (\ref{ZION1}) and~(\ref{ZION2}), we obtain
\begin{align}
\label{ZION3}
{\mathtt P}_{\theta_x\omega}\Big[\dd(Z^{n},
\theta_{n^{-1}\tau_{\mathcal{G}}}M^{n})>\eta\Big]
&\leq {\mathtt P}_{\theta_x\omega}[\tau_{\mathcal{G}}> h_\eps n]
+2{\mathtt P}_{\theta_x\omega}\Big[\sup_{t\in [0, n^{-1}\tau_\G]}|
Z^{n}|>\frac{\eta}{4}\Big]
 \nonumber\\
 &\phantom{**}
+{\mathtt E}_{\theta_x\omega}
\Big(\mathtt{P}_{\theta_{X_{\tau_\G}}\omega}\Big[\sup_{t\in [0,
h_\eps]}|Z^{n}|>\frac{\eta}{4}\Big]
 \nonumber\\
 &\phantom{*****}
+\1{\tau_{\mathcal{G}}\leq h_\eps n}
\mathtt{P}_{\theta_{X_{\tau_\G}}\omega}
\Big[\dd(Z^{n},\theta_{n^{-1}\tau_{\mathcal{G}}}Z^{n})
>\frac{\eta}{4}\Big]\Big).
\end{align}
By definition of 
a $(\eps,n)$-good site, we can choose $\eps>0$ small enough in such a way that 
$\eps\leq \min\{\eta/4, \tilde{\eps}(32\|F\|_{\infty})^{-1}\}$ and $h_\eps \leq \eta^2/16$. Therefore, we have uniformly in $x\in
[-H\sqrt{n},H\sqrt{n}]$,
\begin{align}
\label{XION}
{\mathtt E}_{\theta_x\omega}
\Big(\mathtt{P}_{\theta_{X_{\tau_\G}}\omega}\Big[\sup_{t\in [0,
h_\eps]}|Z^{n}|>\frac{\eta}{4}\Big]
&
+\1{\tau_{\mathcal{G}}\leq
h_\eps n}\mathtt{P}_{\theta_{X_{\tau_\G}}\omega}
\Big[\dd(Z^{n},\theta_{n^{-1}\tau_{\mathcal{G}}}Z^{n})
>\frac{\eta}{4}\Big]\Big)\leq
\frac{\tilde{\eps}}{32\|F\|_{\infty}}
\end{align}
for all sufficiently large~$n$.
On the other hand, by Lemma~\ref{reacgoodsitefrombad}, we have uniformly in $x\in [-H\sqrt{n},H\sqrt{n}]$,
\begin{equation}
\label{XION1}
{\mathtt P}_{\theta_x\omega}[\tau_{\mathcal{G}}> h_\eps n]
\leq \frac{\tilde{\eps}}{32\|F\|_{\infty}}
\end{equation}
and 
\begin{equation}
\label{XION2}
{\mathtt P}_{\theta_x\omega}\Big[\sup_{t\in [0, n^{-1}\tau_\G]}|
 Z^{n}|>\frac{\eta}{4}\Big]\leq
\frac{\tilde{\eps}}{32\|F\|_{\infty}}
\end{equation}
for sufficiently large $n$.
Combining~(\ref{XION}), (\ref{XION1}), (\ref{XION2}) 
with~(\ref{ZION3}), (\ref{ZION}), (\ref{ZION0}) and~(\ref{ZION-1}),
we have
\[
\Big|{\mathtt E}_{\theta_x\omega}\Big(\!F(Z^{n})
- {\mathtt E}_{\theta_{X_{\tau_{\G}}} \omega}[F(Z^{n})]\Big)\Big| \leq \tilde{\eps}/2.
\]
Together with~(\ref{UFin}), 
we obtain that $R\leq \tilde{\eps}$ which concludes the proof
of Proposition~\ref{UCLTprop1}.
\qed

Next, we prove
\begin{prop}
\label{portmantpartial}
The first statement implies the second one:
\begin{itemize}
\item[(i)] for any $F\in \C^u_b(C(\R_+),\R)$, we have $\IP$-a.s., 
\[
\lim_{n \to \infty} \sup_{x\in [-H\sqrt{n},H\sqrt{n}]}\Big| 
 {\mathtt E}_{\theta_x \omega}[F(Z^{n})]-E[F(W)]\Big|=0;
 \]
 \item[(ii)] for any open set $G\subset C(\R_+)$, we have $\IP$-a.s., 
\[
\liminf_{n \to \infty} \inf_{x\in [-H\sqrt{n},H\sqrt{n}]}
{\mathtt P}_{\theta_x \omega}[Z^{n}\in G]\geq P[W\in G].
 \]
\end{itemize}
\end{prop}

\noindent
\textit{Proof.} Let $G$ be an open
set. Then, there exists a sequence $(F_k,k\geq 1)\subset 
\C^u_b(C(\R_+),\R)$ such that $F_k \uparrow {\bf 1}_G$ pointwise as
$k \to \infty$. Thus, we have for all $\omega$, $n$, $k$ and $x\in
[-H\sqrt{n},H\sqrt{n}]$
\begin{equation}
\label{Port11}
{\mathtt P}_{\theta_x \omega}[Z^{n}\in G]\geq 
{\mathtt E}_{\theta_x \omega}[F_k(Z^{n})].
\end{equation}
Then, fix $\eps>0$. By the monotone convergence theorem, 
there exists $k_0$ such that for all $k\geq k_0$,
\begin{equation}
\label{Port22}
E[F_k(W)]\geq P[W\in G]-\frac{\eps}{2}.
\end{equation}
Now, by (i), $\IP$-a.s., for all $k\geq k_0$, we have that 
for $n\geq n_0(k,\omega)$ and all $x\in [-H\sqrt{n},H\sqrt{n}]$,
\begin{equation}
\label{Port33}
{\mathtt E}_{\theta_x \omega}[F_k(Z^{n})]\geq E[F_k(W)]
-\frac{\eps}{2}.
\end{equation}
Combining~(\ref{Port11}) and~(\ref{Port33}), we have, $\IP$-a.s., 
for all $k\geq k_0$, for all $n\geq n_0(k,\omega)$ and all $x\in
[-H\sqrt{n},H\sqrt{n}]$
\begin{equation}
\label{Port44}
{\mathtt P}_{\theta_x \omega}[Z^{n}\in G]\geq E[F_k(W)]
-\frac{\eps}{2}.
\end{equation}
Then, combining~(\ref{Port22}) and~(\ref{Port44}) we obtain $\IP$-a.s.,
for all sufficiently large~$n$,
\begin{equation}
\inf_{x\in [-H\sqrt{n},H\sqrt{n}]}{\mathtt P}_{\theta_x \omega}[Z^{n}\in G]\geq P[W\in G]-\eps.
\end{equation}
As $\eps$ is arbitrary, take the $\liminf_{n\to \infty}$ in the last
inequality to show that (i)~$\Rightarrow$~(ii).
\qed
Then, in the next proposition, we show that ``for every'' and
``$\IP$-a.s.'' can be interchanged:
\begin{prop}
\label{portlast}
The following statements are equivalent:
\begin{itemize}
\item[(i)] we have $\IP$-a.s., for every open set~$G$,
\[
\liminf_{n \to \infty} 
\inf_{x\in [-H\sqrt{n},H\sqrt{n}]}{\mathtt P}_{\theta_x
\omega}[Z^{n}\in G]\geq P[W\in G];
 \]
\item[(ii)] for every open set~$G$, we have $\IP$-a.s., 
\[
\liminf_{n \to \infty} 
\inf_{x\in [-H\sqrt{n},H\sqrt{n}]}{\mathtt P}_{\theta_x
\omega}[Z^{n}\in G]\geq P[W\in G].
 \]
\end{itemize}
\end{prop}

\noindent
\textit{Proof.} 
We only have to
show that (ii)~$\Rightarrow$~(i). Since a basis of open sets for the topology of $C(\R_+)$ is formed by open balls of rational radii about piecewise linear functions connecting rational points, there exists a
countable family $\mathcal{G}$ of open sets such that for every open
set~$G$ there exists a sequence $(O_n,n\geq 1) \subset
\mathcal{G}$ such that ${\bf 1}_{O_n} \uparrow {\bf 1}_{G}$ pointwise
as $n\to \infty$. By~(ii), since the family $\mathcal{G}$ is
countable we have, $\IP$-a.s., for all $O\in \mathcal{G}$, 
\begin{equation}
\liminf_{n \to \infty} \inf_{x\in [-H\sqrt{n},H\sqrt{n}]}
{\mathtt P}_{\theta_x \omega}[Z^{n}\in O]\geq P[W\in O].
\end{equation}
Then, the same kind of reasoning as that used in the proof of
Proposition~\ref{portmantpartial} would
provide the desired result.
\qed

Now we are ready to finish the proof of Theorem~\ref{Theouni}.
Recall that we have to prove that the following statements hold:
\begin{itemize}
\item[(i)] we have $\IP$-a.s., for any $F\in \C_b(C(\R_+),\R)$, 
\[\lim_{n \to \infty} \sup_{x\in [-H\sqrt{n},H\sqrt{n}]}\Big| 
 {\mathtt E}_{\theta_x \omega}[F(Z^{n})]-E[F(W)]\Big|=0;
 \]
\item[(ii)] we have $\IP$-a.s., for any $F\in \C^u_b(C(\R_+),\R)$, 
\[
\lim_{n \to \infty} \sup_{x\in [-H\sqrt{n},H\sqrt{n}]}\Big| 
 {\mathtt E}_{\theta_x \omega}[F(Z^{n})]-E[F(W)]\Big|=0;
 \]
 \item[(iii)] we have $\IP$-a.s., for any closed set $B$,
\[
\limsup_{n \to \infty} \sup_{x\in [-H\sqrt{n},H\sqrt{n}]}
{\mathtt P}_{\theta_x \omega}[Z^{n}\in B]\leq P[W\in B];
 \]
\item[(iv)]  we have $\IP$-a.s., for any open set $G$,
\[
\liminf_{n \to \infty} \inf_{x\in [-H\sqrt{n},H\sqrt{n}]}{\mathtt
P}_{\theta_x \omega}[Z^{n}\in G]\geq P[W\in G];
 \]
 \item[(v)]  we have $\IP$-a.s., for any $A\in \mathcal{B}$ 
such that $P[W\in \partial
A]=0$,
\[
\lim_{n \to \infty} \sup_{x\in [-H\sqrt{n},H\sqrt{n}]}\Big| 
 {\mathtt P}_{\theta_x \omega}[Z^{n}\in A]-P[W\in A]\Big|=0.
 \]
\end{itemize}
Essentially, we follow the proof of
Theorem 2.1 of~\cite{Billing}.  Of course, (i)~$\Rightarrow$~(ii) is
trivial. The proof of the fact that (ii) $\Rightarrow$ (iii) (and, by complementation, that (ii) $\Rightarrow$ (iv))
 is a consequence of Propositions~\ref{UCLTprop1}, \ref{portmantpartial} and \ref{portlast}.
Let us show that (iii)~$\Leftrightarrow$~(v). We start by
showing (iii) $\Rightarrow$ (v). 
Let $\mathring{A}$ denote the interior of~$A$ and~$\bar{A}$ 
denote its closure. If~(iii) holds, then so does~(iv), and hence
$\IP$-a.s.,
\begin{align}
\label{Port5}
P[W\in \bar{A}]
&\geq \limsup_{n \to \infty} 
\sup_{x\in [-H\sqrt{n},H\sqrt{n}]}{\mathtt P}_{\theta_x
\omega}[Z^{n}\in \bar{A}]\geq \limsup_{n \to \infty} \sup_{x\in
[-H\sqrt{n},H\sqrt{n}]}{\mathtt P}_{\theta_x \omega}[Z^{n}\in
A]\nonumber\\
&\geq \liminf_{n \to \infty} 
\sup_{x\in [-H\sqrt{n},H\sqrt{n}]}{\mathtt P}_{\theta_x
\omega}[Z^{n}\in A]\geq \liminf_{n \to \infty} \sup_{x\in
[-H\sqrt{n},H\sqrt{n}]}{\mathtt P}_{\theta_x \omega}[Z^{n}\in
\mathring{A}]\nonumber\\
& \geq P[W\in \mathring{A}].
\end{align}
Since $P[W\in \partial A]=0$, the first and the last terms
in~(\ref{Port5}) are both 
equal to $P[W\in A]$, and~(v) follows.
We continue by showing that (v)~$\Rightarrow$~(iii). Since
$\partial\{ w\in C(\R_+): \dd(w,B)\leq \delta\}=\{w\in C(\R_+):
\dd(w,B)=\delta\}$, the sets $\partial\{ w\in C(\R_+): \dd(w,B)\leq
\delta\}$ are disjoint for distinct $\delta$, hence at most countably
many of them can have positive $P[W\in \cdot~]$-measure. Thus, for
some sequence of positive $\delta_k$ such that $\delta_k\to 0$ as
$k\to \infty$, the sets $B_k=\{w:\dd(w,B)\leq \delta_k\}$ are such
that $P[W\in \partial B_k]=0$. If (v) holds, then we can apply the
same sequence of arguments used to show Proposition \ref{portmantpartial}
with the sequence $({\bf1}_{B_k}, k\geq 1)$ instead of $(F_k, k\geq 1)$. 

Finally, 
we show
that 
(iii)~$\Rightarrow$~(i). 
Suppose that~(iii) holds and that $F\in \C_b(C(\R_+),\R)$. By transforming $F$ linearly (with positive coefficient for the first-degree term) we can reduce the problem to the case in which $0\leq F<1$. For a fixed 
integer~$k$, let $B_i$ be the closed set $B_i=\{w:i/k\leq F(w)\}$,
$i=0,\dots, k$. Since $0\leq F< 1$, we have for all~$\omega$, $n$
and all $x\in [-H\sqrt{n},H\sqrt{n}]$
\[
 {\mathtt E}_{\theta_x \omega}[F(Z^n)]
 <\sum_{i=1}^{k}\frac{i}{k}{\mathtt P}_{\theta_x
\omega}\Big[\frac{i-1}{k}\leq F(Z^n)<\frac{i}{k}\Big]
\]
which implies
\[
 {\mathtt E}_{\theta_x \omega}[F(Z^n)]
 <\frac{1}{k}+\frac{1}{k}\sum_{i=1}^{k}{\mathtt P}_{\theta_x
\omega}[Z^n\in B_i].
\]
If~(iii) holds then we have $\IP$-a.s., $\limsup_{n\to
\infty}\sup_{x\in [-H\sqrt{n},H\sqrt{n}]} {\mathtt P}_{\theta_x
\omega}[Z^n\in B_i]\leq P[W\in B_i]$ for all~$i$, hence we can
deduce that we have $\IP$-a.s., 
\[
\limsup_{n \to \infty}\sup_{x\in [-H\sqrt{n},H\sqrt{n}]} 
{\mathtt E}_{\theta_x \omega}[F(Z^n)]\leq \frac{1}{k}+ E[F(W)].
\]
Since $k$ is arbitrary, we obtain for all $F\in \C_b(C(\R_+),\R)$,
\begin{equation}
\label{Port6}
\limsup_{n \to \infty} \sup_{x\in [-H\sqrt{n},H\sqrt{n}]} {\mathtt
E}_{\theta_x \omega}[F(Z^n)]\leq E[F(W)].
\end{equation}
Note that (\ref{Port6}) implies
\begin{equation}
\label{Port7}
\limsup_{n \to \infty}\inf_{x\in [-H\sqrt{n},H\sqrt{n}]}
{\mathtt E}_{\theta_x \omega}[F(Z^n)]\leq E[F(W)].
\end{equation}
Applying (\ref{Port7}) to $(-F)$ yields 
$\liminf_{n \to \infty}\sup_{x\in [-H\sqrt{n},H\sqrt{n}]} {\mathtt
E}_{\theta_x \omega}[F(Z^n)]\geq E[F(W)]$ 
which together with~(\ref{Port6}) implies (i), and thus the
proof of Theorem~\ref{Theouni} is concluded.
\qed
%


\section*{Acknowledgments}
C.G.\ is grateful to FAPESP (grant 2009/51139--3) for financial
support. S.P.\ was partially supported by
CNPq (grant 300886/2008--0). 
Both also thank CNPq (472431/2009--9) and FAPESP (2009/52379--8)  
for financial support.
The authors thank the referee for very careful reading and lots of suggestions that permitted to considerably improve the paper.

\end{document}